\documentclass[11pt]{article}
\usepackage{times}
\usepackage{indentfirst}
\usepackage{amsthm,amsmath, amssymb}
\usepackage{calrsfs}
\usepackage{graphicx}
\usepackage{multicol}
\usepackage{caption}
\usepackage{subcaption}

\usepackage{algorithm}
\usepackage[noend]{algpseudocode}

\makeatletter
\def\BState{\State\hskip-\ALG@thistlm}
\makeatother

\newtheorem{Theorem}{Theorem}
\newtheorem{Lemma}{Lemma}

\newtheorem{Corollary}{Corollary}
\newtheorem{Remark}{Remark}

\newcommand{\eq}[1]{(\ref{eq:#1})}
\newcommand{\lemt}[1]{Lemma~\ref{lem:#1}}
\newcommand{\thr}[1]{Theorem~\ref{thr:#1}}

\newcommand{\sectn}[1]{Section~\ref{sect:#1}}
\newcommand{\rem}[1]{Remark~\ref{rem:#1}}

	{\par\noindent{\bf Proof.}} 
	{\hfill$\scriptstyle\blacksquare$} 

\newcommand{\PB}{\mathbb{P}}
\newcommand{\Expect}{\mathbb{E}}

\begin{document}

\captionsetup[figure]{labelfont={bf},name={Fig.},labelsep=period}
\title{\bf \Large Stochastic stability of a system of perfect integrate-and-fire inhibitory neurons}

\author{Timofei Prasolov\\Heriot-Watt University}


\maketitle

\begin{abstract}
We study a system of perfect integrate-and-fire inhibitory neurons. It is a system of stochastic processes which interact through receiving an instantaneous increase at the moments they reach certain thresholds. In the absence of interactions, these processes behave as a spectrally positive L\'evy processes. Using the fluid approximation approach, we prove convergence to a stable distribution in total variation.
\end{abstract}

\begin{quotation}
\end{quotation}


\section{Introduction}

We analyse the stochastic stability of a model of a neural network. Our model is inspired by the \textit{stochastic integrate-and-fire neuron} model. The original model of membrane potentials was introduced by Lapicque (1907) and has been developed over the years (for a review of the model see, e.g., Burkitt (2006a,b)). In this model, at any time $t$, the internal state of a neuron $i$ is given by its membrane potential $Z_i(t)$, which evolves according to a stochastic differential equation
$$dZ_i(t) = F(Z_i(t), I(t), t)dt + \sigma(Z_i(t), I(t), t)dW_i(t),$$
where $F$ is a drift function, $\sigma$ the diffusion coefficient, $I$ is the neuronal input, and $W_i$ is a Brownian motion (see, e.g., Gerstner and Kistler (2002)). The process $W_i(t)$ represents combined internal and external noise. The process $I$ models firings of neurons' potentials (or ''spikes''): whenever a potential $Z_i(t)$ reaches certain threshold, it resets to a base-level, and the neuron sends signals to other neurons. 

A large number of experiments have given us an understanding of the dynamics of a single neuron. For example, Hodgkin and Huxley (1952) found three different types of ion current flowing through a neuron's membrane, and introduced a detailed model of a membrane potential. To give a basic description, without any input, the neuron is at rest, corresponding to a constant membrane potential. Given a small change, the membrane potential returns to the resting position. If the membrane potential is given a big enough increase, it reaches a certain threshold, and exhibits a pulse-like excursion that will effect connected neurons. After the pulse, the membrane potential does not directly return to the resting potential, but goes below it. This is connected to the fact that a neuron can not have two spikes one right after another.      

As for the neural networks, neurons form a connected graph with synapses between neurons. When a presynaptic neuron fires a spike it sends a signal through a synapse to a postsynaptic neuron. A neuron is called \textit{inhibitory} if its signals predominantly move the membrane potentials away from a threshold; and \textit{excitatory} if they move potentials toward a threshold. In this paper, we consider a model containing inhibitory neurons only. It is important to point out that the effect of a signal depends on the potential of a receiving neuron. For example, if the membrane potential of a postsynaptic neuron is lower than that of a corresponding inhibitory synapse, the effect of a signal will be reversed. Therefore, in the models where signals and potentials are assumed to be independent, it is important to assume that the potentials should not decrease too much. 

One of the most classical models, introduced by Stein (1965), is \textit{leaky integrate-and-fire neuron} model, where
$$F(Z_i(t), I(t), t) = -\alpha Z_i(t) + I(t), \ \  \sigma(Z_i(t), I(t), t) = \sigma = \text{const.}$$
There are several variations of this model. For instance, nonlinear models were considered, such as the quadratic model (see, e.g., Latham \textit{et al.} (2000)) where $-\alpha Z_i(t)$ is replaced with\\
 $a(Z_i(t) - z_{rest})(Z_i(t) - z_c)$, where $z_c > z_{rest}$. Another direction for generalisation of this model is the Spike Response Model (see, e.g., Gerstner and Kistler (2002), Chapter 4.2). In this model, the relation between the dynamics and the potential is determined by the time of the last spike. This allows one to explicitly forbid spikes to occur one right after another and to write the dynamics in integrated form. 

In this paper we consider \textit{perfect integrate-and-fire neuron} model, where $\alpha = 0$ and, therefore, the decay of the membrane potential over time is neglected. This restriction is a stepping stone to achieve more general results and it allows us to write the model in integrated form. 

In our model, the spikes and corresponding signals are represented by shifts from a threshold of a random length, independent of everything else. We analyse the system under certain conditions on the distribution of those shifts and prove stability. Instead of considering the recurrence of sets $[-k, H]^N$ (where $H$ is a threshold and $N$ is a number of neurons), we move each coordinate down and reflect the system to work with more convenient sets $[0, k+H]^N$. Thus, in our model, membrane potentials are nonnegative processes that jump to a random positive level after reaching zero. Signals from inhibitory neurons push membrane potentials from the threshold, i.e. they are positive shifts. It is important to note that we assume that travel time of signals between neurons is zero, which in general can cause uncertainty in the order of spikes. However the inhibitory signals do not cause spikes right away, and we assume that the potentials $Z_i(t)$ almost surely do not reach their thresholds at the same time. We refer to Taillefumier \textit{et al.} (2012) for further discussion.

It is often assumed that the studied system of neurons is itself a part of a much larger system of neurons. The corresponding effect on our system is often modelled by a multivariate Brownian motion $W(t)$ with a drift (the drift guaranties the stability of a system of a single neuron). However, we can generalise it to a multivariate spectrally positive (i.e. with positive jumps) L\'evy process $X(t)$ to account for inhibitory signals. It is important for our analysis that the signals do not influence the dynamics of the process $Z(t)$ if it is away from the threshold, i.e. we have $dZ(t) = dX(t)$ if $Z_i(t) > 0$, for $i\in \{1, \ldots, N\}$. Nevertheless, the number of spikes $\eta_i(t)$ before time $t$ is essential to stability analysis. The fact that $\eta_i(t)$ is not pathwise monotone with respect to signal sizes or the initial state brings certain difficulties in proving stability. 

As was mentioned, a system of a single neuron is stable, however, for a general distribution of signals between neurons, ''partial stability'' can occur when only a strict subset of neurons (maybe random) stabilises, while membrane potentials of other neurons are ''pushed'' to infinity (which contradicts the physical setup). The latter is of independent mathematical interest and will not be discussed in the main body of the paper.

Under specific conditions on average signals and the drift $\Expect X(1)$, we prove the positive recurrence of the system using the fluid approximation approach, introduced by Rybko and Stolyar (1992) and Dai (1995) (see also Stolyar (1995)). Although this method is usually shown on queueing networks, it is quite universal, and applicable to our model too. Using results from Section 7 of Borovkov and Foss (1992) (see, also, Chapter VII of Asmussen (2003)),  we prove \textit{Harris positive recurrence} and convergence to stationary distribution in total variation. We refer to Foss and Konstantopoulos (2004) for an overview of some stochastic stability methods.

The \textit{stochastic integrate-and-fire neuron} model has received an increasing amount of attention in recent years. There are a number of papers considering mean-field limits of such systems. De Masi \textit{et al.} (2015) consider a model with identical inhibitory neurons, where each membrane potential has a drift to the average potential. Inglis and Talay (2015) consider general signals between neurons and describe signal transmissions through the use of the cable equation (instead of instant transmissions). Robert and Touboul (2016) consider a model where neurons do not have a fixed threshold and spikes occur as a inhomogeneous Poisson process, with intensity given as a function of a membrane potential, and prove ergodicity. Several authors studied cases when excitatory signals lead to so-called \textit{blow-up phenomena} (e.g. C\'aceres \textit{et al.}(2011), Delarue \textit{et al.} (2015)).

This paper is structured as follows. In \sectn{model&results}  we define our model, introduce auxiliary concepts and notations, and formulate our results. In particular, in \sectn{fluid_model} we introduce the fluid model and formulate related technical results. In \sectn{auxiliary_results} we prove important auxiliary results. In \sectn{lemma_positive_recurrence} we prove positive recurrence. In \sectn{lemma_mixing} we prove that our model satisfies the classical ''minorization'' condition. The Appendix includes the remaining auxiliary results and comments.

\section{Model and results}\label{sect:model&results}

We analyse a network of $N$ stochastic perfect integrate-and-fire inhibitory neurons. At any time $t$, the internal state of all neuron is given by a multidimensional process $Z(t)$ which represents neurons' membrane potential. Let $X(t)$ be a $N$-dimensional spectrally positive left-continuous L\'evy process with a finite mean and its distribution has a non-degenerate absolute continuous component. The process $X(t)$ represents combined internal and external noise. Let $\nu_i = - \Expect X_i(1) > 0$ and $X^0_i(t) = \nu_i t + X_i(t)$. While $Z(t) \in (0, \infty)^N$, membrane potentials evolve as the process $X(t)$, i.e. $dZ(t) = dX(t)$. Let $\{\{\xi^{(k)}_{ij}\}_{i, j=1}^N\}_{k=1}^\infty$ be i.i.d. random matrices, independent of everything else, with a.s. strictly positive elements.

\begin{Remark}
One can allow the absolute continuous component of the distribution of the process $X_i(t)$ to be degenerate (for example, take a sum of a Poisson process and a linear function $-at$) and, instead, condition the distribution of the matrix $\{\xi^{(1)}_{ij}\}_{i, j=1}^N$ to have an absolute continuous component. The main result of the paper would still hold and the proof would need few minor changes.
\end{Remark}

Let $b_{ij} = \Expect \xi^{(1)}_{ij} < \infty$ and $S_{ij}(n) = \sum_{k=1}^{n} \xi^{(k)}_{ij},$ for $i, j \in \{1, 2, \ldots, N\}$. If potential $Z_i(t)$ hits non-positive values for the $k$-th time, then instantaneously it increases to $\xi^{(k)}_{ii}$ and other membrane potentials increase by $\xi^{(k)}_{ij}$. We call this event ''a spike of neuron $i$``.

Let $Z^{\textbf{z}}(t) = (Z^{\textbf{z}}_1(t), \ldots, Z^{\textbf{z}}_N(t)) \in \mathcal{Z} = [0, \infty)^N$ be the membrane potentials at time $t$ with an initial value $\textbf{z} = (z_1, \ldots, z_N)$. Let $T^{\textbf{z}}_{i0} = 0$ and
$$T^{\textbf{z}}_{ik} = \inf\{t> T^{\textbf{z}}_{i(k-1)}: \ Z^{\textbf{z}}_i(t) \leq 0\}, \ \ \text{for $k\geq 1$,}$$ 
be the times when neuron $i$ reaches its threshold. Let $\eta^{\textbf{z}}_i(0) = 0$ and $\eta^{\textbf{z}}_i(t) = \max\{k: \ T^{\textbf{z}}_{ik} < t\}$ be the number of spikes of $Z^{\textbf{z}}_i(t)$ before time $t$.  Then the dynamics of the system is given by
\begin{equation}\label{eq:dynamic_equation}
Z^{\textbf{z}}_i(t) = z_i + X_i(t) + \sum_{j=1}^N \sum_{k=1}^{\eta^{\textbf{z}}_j(t)} \xi^{(k)}_{ji} = z_i + X_i(t) + \sum_{j=1}^N S_{ji}(\eta^{\textbf{z}}_j(t)), \ \ i=1, \ldots, N.
\end{equation}

Before talking about stability of the system it is important to point out that, due to the negative drift, one can easily show that potential of an isolated neuron is stable (this is also a subcase of our main result).

\begin{Remark}\label{rem:partial_stability}
There are examples of parameters $\nu_i$ and $b_{ij}$ such that there exists a subset of neurons which, after reaching stability, can ''push other neurons to infinity``. We do not discuss such cases of partial stability in the main body of the paper, however, we include a few comments in the Appendix. 
\end{Remark}

We assume that all potentials have the same drift $\nu$ and that signals from neuron $i = 1, \ldots, N$ to all other neurons have the same mean $w_i$. We have
\begin{equation}\label{eq:simple_network}
\nu_i = \nu >0, \ \ b_{ij} = \Expect \xi^{(1)}_{ij} = w_i > 0 \ \ \text{and} \ \ b_{ii} = H_i > w_i, \ \ \text{for $i=1, \ldots, N$ and $j\neq i$.}
\end{equation}

\begin{Theorem}\label{thr:ergodicity}
Assume condition \textnormal{\eq{simple_network}} to hold. Then the process  $(Z^{\textbf{z}}_1(t), \ldots, Z^{\textbf{z}}_N(t))$ is Harris positive recurrent: there is a distribution $\pi$ such that
$$\sup_{A}|\PB\{Z^{\textbf{z}}(t) \in A\} - \pi(A)| \to 0 , \ \text{as $t\to\infty$.}$$
\end{Theorem}

\begin{Remark}\label{rem:spike_rate}
Given \textnormal{\eq{simple_network}}, matrix $B = (b_{ij})_{i, j=1}^N$ is invertible and
$$(\textbf{1}B^{-1})_i = \frac{1}{\left( H_i - w_i \right) \left( 1 + \sum_{k=1}^{N}\frac{w_k}{H_k - w_k} \right)}, \ \ \text{for $i = 1, \ldots, N$,}$$
where $\textbf{1} = (1, \ldots, 1)$. Vector $\nu \textbf{1}B^{-1}$ represents rates of spikes when stability is achieved. In particular, for large $t > \nu^{-1}(1 + \sum_{k=1}^N w_k / (H_k - w_k))$ and for each sequence $\textbf{z}_n$, $\|\textbf{z}_n\| \to \infty$, there exists a subsequence $\textbf{z}_{n_k}$ such that
$$\frac{\eta^{\textbf{z}_{n_k}}(\|\textbf{z}_{n_k}\|(t+\Delta)) - \eta^{\textbf{z}_{n_k}}(\|\textbf{z}_{n_k}\|t)}{\|\textbf{z}_{n_k}\|\Delta} \Rightarrow \nu \textbf{1}B^{-1}, \ \text{for $\Delta >0$.}$$  
\end{Remark}

We prove \thr{ergodicity} following two standard steps. For the reader's convenience, we formulate those steps as lemmas. Let $\tau^{\textbf{z}}(\varepsilon, B) = \inf\{t> \varepsilon: \ Z^{\textbf{z}}(t) \in B\}$ be the first hitting time of a set $B$ after time $\varepsilon$. The first step is the proof of positive recurrence which we achieve via the fluid approximation method.
\begin{Lemma}\label{lem:positive_recurrence}
There exists $k_0 >0$ such that for $V = \{\textbf{z}\in \mathcal{Z}: \|\textbf{z}\| < k_0\}$ we have
$$\sup_{\textbf{z}\in V} \Expect\tau^{\textbf{z}}(\varepsilon, V) < \infty.$$
\end{Lemma}
In the second step, we show that our model satisfies the classical ''minorization'' condition.
\begin{Lemma}\label{lem:minorant_measure}
There exist a number $p>0$ and a probability measure $\psi$ such that for a uniformly distributed r.v. $U \in [1,2]$, independent of everything else, we have 
$$\inf_{\textbf{z} \in V} \PB\{Z^{\textbf{z}}(U) \in B\} \geq p \psi(B).$$
\end{Lemma}
Using Lemmas $1$ and $2$ we can prove that conditions of Theorem 7.3 from Borovkov and Foss (1992) are satisfied, which gives us the result. The proof of \lemt{positive_recurrence} is based on the fluid approximation. We dedicate the following subsection to formulate corresponding definitions and auxiliary results. We point out that we need to assume condition \eq{simple_network} only in the proof of \lemt{positive_recurrence} and in \rem{spike_rate}. 

One of the difficulties of our model is lack of path-wise monotonicity for the number of spikes $\eta^{\textbf{z}}(t)$ with respect to signals $\xi^{(k)}_{ij}$ or initial state $\textbf{z}$. In general, making one neuron firing a spike earlier may lead to other spikes occur later. However, there is a ''partial monotonicity'' which allows us to get an upper bound for process $\eta^{\textbf{z}}(t)$ with useful properties.

Since all neurons are inhibitory, one way to increase the number of spikes is to remove all interactions between neurons. Let the process $\widetilde{Z}^{\textbf{z}}$ be the transformation of the process $Z^{\textbf{z}}$ by replacing signals $\xi^{(k)}_{ji}$, $j\neq i$,  by $0$ for $k\geq 1$ (trajectories of $X(t)$ remain the same). The resulting process has a simpler dependence between coordinates and it has a greater number of spikes before any time $t>0$ than that of $Z^{\textbf{z}}$. For our convenience, we want to remove the dependence of the upper bound on $\textbf{z}$ (which is significant because we take $\textbf{z}$ large in the following lemmas) and make the time until the first spike to have the same distribution with the rest of waiting times. Let the process $\bar{Z}$  be the transformation of the process $\widetilde{Z}^{\textbf{z}}$, so that $\bar{Z}_i(0) \overset{d}{=} \xi^{(1)}_{ii}$, $1\leq i \leq N$. Let $\widetilde{\eta}^{\textbf{z}}$ and $\bar{\eta}$ be the number of spikes in processes $\widetilde{Z}^{\textbf{z}}$ and $\bar{Z}$, respectively.

\begin{Lemma}\label{lem:spikes_bound}
We have 
$$\eta^{\textbf{z}}_i(t) \leq \widetilde{\eta}^{\textbf{z}}_i(t), \ \ \text{a.s.,}$$
$$\widetilde{\eta}^{\textbf{z}}_i(t) \overset{st.}{\leq} 1 + \bar{\eta}_i(t),$$ 
and $\bar{\eta}_i(t)$ is an undelayed renewal process, which satisfies the integral renewal theorem and SLLN
$$\frac{\Expect \bar{\eta}_i(t)}{t} \to \frac{\nu_i}{b_{ii}} \overset{a.s.}{\leftarrow} \frac{\bar{\eta}_i(t)}{t}$$
\end{Lemma}

\subsection{Fluid model and corresponding auxiliary results}\label{sect:fluid_model}

Let us define the fluid approximation model. Let $\rho(\textbf{x}, \textbf{y}) = \sum_{i=1}^N |x_i - y_i|$ be the metric on our space $\mathcal{Z}$ and $\|\textbf{x}\| = \rho(\textbf{x}, 0)$, for $\textbf{x}, \textbf{y} \in \mathcal{Z}$. For each $\textbf{z}\in\mathcal{Z}$, introduce a family of scaled processes
$$\widehat{Z}^{\textbf{z}} = \left\lbrace \widehat{Z}^{\textbf{z}}(t) = \frac{Z^{\textbf{z}}(\|\textbf{z}\|t)}{\|\textbf{z}\|}, \ t\geq 0\right\rbrace .$$ We call the family
$$\widehat{Z} = \{\widehat{Z}^\textbf{z}, \ \|\textbf{z}\| \geq 1\}$$
\textit{relatively compact (at infinity)} if, for each sequence $\widehat{Z}^{\textbf{z}_n}$, $\|\textbf{z}_n\| \to \infty$, there exists a subsequence $\widehat{Z}^{\textbf{z}_{n_k}}$ that converges weakly (in Skorokhod topology) to some limit process $\varphi^Z = \{\varphi^Z(t), \ t\geq 0\}$, which is called a \textit{fluid limit}. A family of such limits is called a \textit{fluid model}. The fluid model is \textit{stable} if there exists a finite constant $T$ such that $\|\varphi^Z(T)\| = 0$ a.s. for any fluid limit $\varphi^Z$ (there are several equivalent definitions of stability of a fluid model, see e.g. Stolyar (1995)). Based on stability of a fluid model, one can prove positive recurrence of the original Markov process following the lines of Dai (1995).

Using \lemt{spikes_bound} we prove the next result.

\begin{Lemma}\label{lem:relative_compact}
The family of processes $\{Z^{\textbf{z}}, \ \textbf{z}\in\mathcal{Z}\}$ is such that
\begin{itemize}
\item for all $t>0$ and $\textbf{z}\in\mathcal{Z}$,
$$\Expect\|Z^\textbf{z} (t)\| < \infty$$
and moreover, for any K,
$$\sup_{\|\textbf{z}\|\leq K}\Expect\|Z^\textbf{z} (t)\| < \infty;$$
\item for all $0\leq u < t$, the family of random variables 
$$\{\rho(\widehat{Z}^\textbf{z}(u), \widehat{Z}^\textbf{z}(t)), \ \|\textbf{z}\|\geq 1\}$$
is uniformly integrable and there exists a constant $C$ such that
$$\limsup_{\|\textbf{z}\|\to\infty} \PB\{\sup_{u',t'\in [u, t]}\rho(\widehat{Z}^\textbf{z}(u'), \widehat{Z}^\textbf{z}(t')) > C(t-u) \} = 0.$$
\end{itemize}
\end{Lemma}

With this result, one can follow the lines of the proof of Theorem 7.1 from Stolyar (1995) to obtain the following.

\begin{Corollary}
The family of processes $\widehat{Z}$ is relatively compact and every fluid limit $\varphi^Z$ is an a.s. Lipschitz continuous function with Lipschitz constant $C+1$.
\end{Corollary}
Additionally, the function $\varphi^Z(t)$ is a.s. differentiable. We call a time $t_0$ a \textit{regular point} if $\varphi^Z(t)$ is differentiable at $t_0$. Further more, we have
\begin{equation}\label{eq:lipschitz_integration}
\varphi^Z(t) - \varphi^Z(s) = \int_{s}^t \frac{d \varphi^Z}{du}(u) du, \ \ t > s > 0,
\end{equation}
where the derivative is arbitrarily defined (for example, it equals zero) outside regular points.

Let $\widehat{\eta}^{\textbf{z}}(t) = \eta^{\textbf{z}}(\|\textbf{z}\|t) / \|\textbf{z}\|$. Following the lines of the proof of \lemt{relative_compact}, one can prove similar results for the family $\widehat{\eta} = \{\widehat{\eta}^\textbf{z}, \ \|\textbf{z}\| \geq 1\}$. Denote a fluid limit of  $\widehat{\eta}$ as $\varphi^{\eta}$. If at time $t$ we have $\varphi^{\eta}_i(t) >0$, then for certain sequence $\textbf{z}_n$ the number of spikes $ \eta^{\textbf{z}_n}_i(\|\textbf{z}_n\|t)$ becomes large. If additionally, $ \eta^{\textbf{z}_n}_i(\|\textbf{z}_n\|t)$ converges to infinity a.s., then by the law of large numbers
$$\frac{S_{ij}(\eta^{\textbf{z}_n}_i(\|\textbf{z}_n\|t))}{\|\textbf{z}_n\|} = \frac{\eta^{\textbf{z}_n}_i(\|\textbf{z}_n\|t)}{\|\textbf{z}_n\|}\frac{S_{ij}(\eta^{\textbf{z}_n}_i(\|\textbf{z}_n\|t))}{\eta^{\textbf{z}_n}_i(\|\textbf{z}_n\|t)} \Rightarrow \varphi^\eta_i(t) b_{ij}, \ \text{as $n\to\infty$.}$$
If $\varphi^\eta_i(t) = 0$, then the number of spikes is not as large and, if we prove that the left-hand side of the last equation converges to zero, the resulting convergence will be of the same form. 

Using this idea we get the following result. 

\begin{Lemma}\label{lem:weak_convergence}
Let $\widehat{\eta}^{\textbf{z}_n}$ converge weakly to a fluid limit $\varphi^{\eta}$ for a sequence $\textbf{z}_n$, $\|\textbf{z}_n\| \to \infty$ as $n\to\infty$. Then we have weak convergence of processes 
$$\left(\left(\frac{1}{\|\textbf{z}_n\|}\sum_{i=1}^N S_{ij}(\eta^{\textbf{z}_n}_i(\|\textbf{z}_n\|t))\right)_{j=1}^N, \ t\geq 0\right) \overset{D}{\Rightarrow} (\varphi^{\eta}(t)B, \ t\geq 0).$$
\end{Lemma}

\section{Proofs of auxiliary results}\label{sect:auxiliary_results}
In this section we prove our auxiliary results for a general matrix $B$ and parameters $\nu_i$.

\subsection{Proof of \lemt{spikes_bound}}

We prove that $T^{\textbf{z}}_{ik} \geq \widetilde{T}^{\textbf{z}}_{ik} $:
\begin{eqnarray*}
T^{\textbf{z}}_{ik} & =&\inf\{t> T^{\textbf{z}}_{i(k-1)}: \ Z^{\textbf{z}}_i(t) \leq 0\} = \inf\{t> T^{\textbf{z}}_{i(k-1)}: \ Z^{\textbf{z}}_i(t) = 0\}\\
& = & \inf\{t> T^{\textbf{z}}_{i(k-1)}: \ z_i + X_i(t) + \sum_{j=1}^N S_{ji}(\eta^{\textbf{z}}_j(t)) = 0\}\\
& = & \inf\{t> T^{\textbf{z}}_{i(k-1)}: \ z_i + X_i(t) + S_{ii}(k-1) + \sum_{j\neq i} S_{ji}(\eta^{\textbf{z}}_j(t)) = 0\}\\
& \geq & \inf\{t> T^{\textbf{z}}_{i(k-1)}: \ z_i + X_i(t) + S_{ii}(k-1) = 0\}.
\end{eqnarray*}

Since $T^{\textbf{z}}_{i0} = \widetilde{T}^{\textbf{z}}_{i0} = 0$, by induction we have
$$T^{\textbf{z}}_{ik} \geq \inf\{t> \widetilde{T}^{\textbf{z}}_{i(k-1)}: \ z_i + X_i(t) + S_{ii}(k-1) = 0\} =  \widetilde{T}^{\textbf{z}}_{ik}.$$
Thus, we get $\eta^{\textbf{z}}_i(t) \leq \widetilde{\eta}^{\textbf{z}}_i(t)$. Since $\widetilde{\eta}^{\textbf{z}}_i(t) - 1$ has the same distribution with $\bar{\eta}_i(t-\widetilde{T}^{\textbf{z}}_{i1}) \leq \bar{\eta}_i(t)$, we have the second inequality. 

The process $\bar{\eta}_i(t)$ is an undelayed renewal process with waiting times having the same distribution with $\tau_i = \inf\{t> 0: X_i (t) = -\xi^{(1)}_{ii}\}$. Using the strong law of large numbers, one can prove that $\Expect \tau_i = b_{ii}/\nu_i$ (see also Borovkov (1965) for a detailed proof). Therefore, via the standard argument of renewal theory the rest of the proof follows (see e.g. Feller (1971)).

\subsection{Proof of \lemt{relative_compact}}
\textbf{Part 1.} Using \lemt{spikes_bound} and positivity of $\xi^{(k)}_{ij}$, we get
$$
\|Z^\textbf{z} (t)\| = \sum_{i=1}^N |z_i + X_i(t) + \sum_{j=1}^N  S_{ji}(\eta^{\textbf{z} }_j(t))| \leq \|\textbf{z}\| + \|X(t)\|
 +\sum_{i=1}^N\sum_{j=1}^N S_{ji}(\widetilde{\eta}^{\textbf{z} }_j(t)).
$$
We have
$$
\{\widetilde{\eta}^{\textbf{z} }_i (t) = m\} = \{-\sum_{k=1}^m \xi^{(k)}_{ii} < z_i + \inf_{0\leq s\leq t} X(s) \leq -\sum_{k=1}^{m-1} \xi^{(k)}_{ii}\}
$$
and, therefore,
$$\widetilde{\eta}^{\textbf{z} }_i (t) = \inf\{m \in \mathbb{Z}^+ : \ \sum_{k=1}^m \xi^{(k)}_{ii} > -z_i - \inf_{0\leq s\leq t} X(s)\}.$$
Since $\{\{\xi^{(k)}_{ij}\}_{i, j=1}^N\}_{k=1}^\infty$ and $(X(t), \ t\geq 0)$ are independent, the random variable $\widetilde{\eta}^{\textbf{z} }_i (t)$ is a stopping time for the sequence $\{\{\xi^{(k)}_{ij}\}_{i, j=1}^N\}_{k=1}^\infty$. By Wald's identity,
$$\Expect S_{ji}(\widetilde{\eta}^{\textbf{z} }_j(t)) = \Expect\sum_{k=1}^{\widetilde{\eta}^{\textbf{z} }_j(t)} \xi^{(k)}_{ji} = \Expect\widetilde{\eta}^{\textbf{z} }_j (t)b_{ji} < \infty.$$

\textbf{Part 2.} We have
\begin{equation}\label{eq:distance_t-u}
\begin{split}
\rho(\widehat{Z}^\textbf{z}(u), \widehat{Z}^\textbf{z}(t)) & = \sum_{i=1}^N \frac{|Z^{\textbf{z}}_i(\|\textbf{z}\|t)-Z^{\textbf{z}}_i(\|\textbf{z}\|u)|}{\|\textbf{z}\|} \leq (t-u)\sum_{i=1}^N \nu_i\\
& + \sum_{i=1}^N \frac{|X^0_i(\|\textbf{z}\|t)-X^0_i(\|\textbf{z}\|u)|}{\|\textbf{z}\|}
 +  \sum_{i=1}^N\sum_{j=1}^N  \frac{S_{ji}(\eta^{\textbf{z}}_j(\|\textbf{z}\|t)) - S_{ji}(\eta^{\textbf{z}}_j(\|\textbf{z}\|u))}{\|\textbf{z}\|}.
\end{split}
\end{equation}
Process $X_i$ is a L\'evy process, from which we have
$$\Expect \frac{|X^0_i(\|\textbf{z}\|t)|}{\|\textbf{z}\|} \leq 2\sup_{0\leq s \leq t}\Expect|X^0_i(s)|, \ \text{for $\|\textbf{z}\| \geq 1$,}$$
and, therefore, the second summand in the right-hand side of \eq{distance_t-u}  is uniformly integrable. By \lemt{spikes_bound}, we have 
$$S_{ij}(\eta^{\textbf{z}}_j(\|\textbf{z}\|t)) - S_{ij}(\eta^{\textbf{z}}_j(\|\textbf{z}\|u)) \overset{st.}{\leq} S_{ij}(1 + \bar{\eta}_i(\|\textbf{z}\|(t-u))).$$
Since $S_{ij}(n)/n \to b_{ij}$ and $\bar{\eta}_i(\|\textbf{z}\|(t-u)) \to \infty$ a.s., we have 
$$\frac{S_{ij}(1 + \bar{\eta}_i(\|\textbf{z}\|(t-u)))}{1 +\bar{\eta}_i(\|\textbf{z}\|(t-u))} \overset{a.s.}{\to } b_{ij},$$
and therefore
\begin{equation*}
\begin{split}
0 & \leq  \frac{S_{ji}(\eta^{\textbf{z}}_j(\|\textbf{z}\|t)) - S_{ji}(\eta^{\textbf{z}}_j(\|\textbf{z}\|u))}{\|\textbf{z}\|} \overset{st.}{\leq} \frac{S_{ij}(1 + \bar{\eta}_i(\|\textbf{z}\|(t-u)))}{\|\textbf{z}\|} \\
& = \frac{1 +\bar{\eta}_i(\|\textbf{z}\|(t-u))}{\|\textbf{z}\|}\frac{S_{ij}(1 + \bar{\eta}_i(\|\textbf{z}\|(t-u)))}{1 +\bar{\eta}_i(\|\textbf{z}\|(t-u))} \to (t-u)\frac{\nu_i}{b_{ii}}b_{ij}
\end{split}
\end{equation*}
a.s. and in $L_1$, as $\|\textbf{z}\|\to\infty$.  Then the distance on the left-hand side of \eq{distance_t-u} is bounded above by the sum of uniformly integrable random variables and therefore is also uniformly integrable.

Given 
$$C>\sum_{i=1}^N\nu_i\left(  1 + \sum_{j=1}^N \frac{b_{ij}}{b_{ii}}\right) ,$$ there exists $\varepsilon >0$ such that for $\|\textbf{z}\|$ large

\begin{equation*}
\begin{split}
\PB\{\sup_{u',t'\in [u, t]}\rho(\widehat{Z}^\textbf{z}(u'), \widehat{Z}^\textbf{z}(t')) > C(t-u) \} & \leq \PB\left\{\sup_{u',t'\in [u, t]}\sum_{i=1}^N \frac{|X^0_i(\|\textbf{z}\|t')-X^0_i(\|\textbf{z}\|u')|}{\|\textbf{z}\|} > \varepsilon(t-u) \right\}\\
 & \leq 2N \PB\left\{ \sup_{s\in [0, t-u]}\frac{|X^0_1(\|\textbf{z}\|s)|}{\|\textbf{z}\|} > \frac{\varepsilon}{2N}(t-u) \right\}  \to 0,
\end{split}
\end{equation*}
by Theorem 36.8 from Sato (1999). 

\subsection{Proof of \lemt{weak_convergence}}
By Skorokhod (1956), it is sufficient to prove that there is a convergence of finite-dimensional distributions on everywhere dense set of times $t$ and that a tightness condition holds. Tightness can be deduced from the second statement of \lemt{relative_compact}. We prove that
\begin{equation}\label{eq:finite_dimensional_convergence}
\PB\left\{\bigcap_{k=1}^K \bigcap_{i, j =1}^N \left\{\frac{S_{ij}(\eta^{\textbf{z}_n}_i(\|\textbf{z}_n\|t_k))}{\|\textbf{z}_n\| } < y^k_{ij}\right\}\right\} \to \PB\left\{\bigcap_{k=1}^K \bigcap_{i =1}^N \left\{\varphi^{\eta}_i(t_k) < \min_{1\leq j\leq N}\frac{y^k_{ij}}{b_{ij}}\right\}\right\}
\end{equation}
as $n\to\infty$, for appropriate $t\geq 0$ and $\textbf{y}\in (0, \infty)^{K N^2}$. 

Define sets
$$C_{ij}^k(n) = \left\{ \frac{S_{ij}(\eta^{\textbf{z}_n}_i(\|\textbf{z}_n\|t_k))}{\|\textbf{z}_n\| } < y^k_{ij} \right\},$$
$$D^k_i (n, m) = \{\eta^{\textbf{z}_n}_i(\|\textbf{z}_n\|t_k) > m\},$$
$$E^k_{ij} (n, \delta) = \left\{\left|  \frac{S_{ij}(\eta^{\textbf{z}_n}_i(\|\textbf{z}_n\|t_k))}{\eta^{\textbf{z}_n}_i(\|\textbf{z}_n\|t_k)} - b_{ij}  \right| \leq \delta  \right\},$$
$$F^{k \pm}_{i} (n, \delta) = \left\{\frac{\eta^{\textbf{z}_n}_i(\|\textbf{z}_n\|t_k)}{\|\textbf{z}_n\|} < \min_{1 \leq j \leq N} \frac{y^k_{ij}}{b_{ij} \mp \delta} \right\}, $$
where $\delta \in (0, \min_{i,j} b_{ij})$. We prove that
$$\PB\left\{F^{k -}_{i} (n, \delta)\right\} + o(1) \leq \PB\left\{\bigcap_{j=1}^N C_{ij}^k(n)\right\} \leq  \PB\left\{F^{k +}_{i} (n, \delta)\right\} + o(1), \ \text{as $n\to \infty$.}$$
For any $\textbf{y}\in (\mathbb{R}^+)^{K N^2}$ such that $(\min_j (y^k_{ij}/b_{ij}))_{i=1}^N$ is a continuity point of the cdf of $(\varphi^{\eta}(t_k))_{k=1}^k$, there is a neighbourhood $\Delta$ of $\textbf{y}$ such that every point $\textbf{x} \in \Delta$ is also a continuity point. Thus, for $\delta$ small we have

\begin{equation*}
\begin{split}
\PB\left\{\bigcap_{k=1}^K \bigcap_{i=1}^N F^{k \pm}_{i} (n, \delta)\right\} & = \PB\left\{\bigcap_{k=1}^K\bigcap_{i =1}^N\left\{\frac{\eta^{\textbf{z}_n}_i(\|\textbf{z}_n\|t_k)}{\|\textbf{z}_n\| } < \min_{1 \leq j \leq N} \frac{y^k_{ij}}{b_{ij} \mp \delta}\right\}\right\} \\
 & \to \PB\left\{\bigcap_{k=1}^K\bigcap_{i =1}^N \left\{\varphi^{\eta}_i(t_k) < \min_{1\leq j\leq N}\frac{y^k_{ij}}{b_{ij} \mp \delta}\right\}\right\}, \ \text{as $n \to \infty$,}
\end{split}
\end{equation*}
and, therefore, by letting $\delta$ converge to $0$, we get \eq{finite_dimensional_convergence}.

By the law of large numbers, we have
$$\PB\{D^k_i (n, m) \cap \overline{E^k_{ij} (n, \delta)}\} \to 0, \ \ \text{as $m\to\infty$,}$$
and
$$\PB\{\overline{C_{ij}^k(n)} \cap \overline{D^k_i (n, m)}\} \to 0, \ \ \text{as $n\to\infty$,}$$
if $m=o(\|\textbf{z}_n\|)$. Take $m=\sqrt{\|\textbf{z}_n\| }$.

From the definitions we have
\begin{multline*}
\left( \bigcap_{j=1}^N C_{ij}^k(n) \cap D^k_i (n, m) \cap E^k_{ij} (n, \delta)\right)  \subseteq \left( F^{k +}_{i} (n, \delta) \cap D^k_i (n, m) \cap E^k_{ij} (n, \delta)\right)\\
 = \left( F^{k +}_{i} (n, \delta) \cap D^k_i (n, m)\right) \setminus \left(  F^{k +}_{i} (n, \delta) \cap D^k_i (n, m) \cap \overline{E^k_{ij} (n, \delta)}\right)  
\end{multline*}
and
$$ \left( F^{k +}_{i} (n, \delta) \cap D^k_i (n, m)\right) = F^{k +}_{i} (n, \delta) \setminus \left( F^{k +}_{i} (n, \delta) \cap \overline{D^k_i (n, m)}\right).$$

Since $m=o(\|\textbf{z}_n\|)$, we have $F^{k +}_{i} (n, \delta) \cap \overline{D^k_i (n, m)} = \overline{D^k_i (n, m)}$ for $n$ large. Combining altogether, we get

\begin{equation*}
\begin{split}
\PB\left\{\bigcap_{j=1}^N C_{ij}^k(n)\right\} & = \PB\left\{\bigcap_{j=1}^N C_{ij}^k(n) \cap  D^k_i (n, m)\right\} + 
	\PB\left\{\bigcap_{j=1}^N C_{ij}^k(n) \cap \overline{D^k_i (n, m)}\right\} \\
& \leq \PB\left\{F^{k +}_{i} (n, \delta)\right\} - \PB\left\{\overline{D^k_i (n, m)}\right\} + \PB\left\{\bigcap_{j=1}^N C_{ij}^k(n) \cap \overline{D^k_i (n, m)}\right\} + o(1),
\end{split}
\end{equation*}
as $n\to\infty$, and
\begin{multline*}
\PB\left\{\bigcap_{j=1}^N C_{ij}^k(n) \cap \overline{D^k_i (n, m)}\right\} -  \PB\left\{\overline{D^k_i (n, m)}\right\} =  \PB\left\{\bigcup_{j=1}^N \overline{C_{ij}^k(n)} \cap \overline{D^k_i (n, m)}\right\} \to 0,
\end{multline*}
as $n\to\infty$. Following the same lines with replacing a set $F^{k+}_i(n, \delta)$  with a set $F^{k-}_i(n, \delta)$ and relations $\subseteq$ and $\leq$ with relations $\supseteq$ and $\geq$, we get the lower bound.

\section{Proof of \lemt{positive_recurrence}}\label{sect:lemma_positive_recurrence}

We prove that under condition \eq{simple_network} fluid limits $\varphi^Z(t)$ are deterministic and uniquely defined by initial value $\varphi^{Z}(0)$. Further, each coordinate of a fluid limit is a continuous piecewise linear function which tends to zero and then remains there.

Let sequence $\textbf{z}_n$, $\|\textbf{z}_n\| \to \infty$, be such that
$$\widehat{Z}^{\textbf{z}_n} \overset{\mathcal{D}}{\Rightarrow} \varphi^Z \ \text{and} \ \widehat{\eta}^{\textbf{z}_n} \overset{\mathcal{D}}{\Rightarrow} \varphi^\eta.$$
By Corollary $1$, function $\varphi^Z$ is a.s. Lipschitz continuous.

Following the lines of the proof of \lemt{weak_convergence}, one can easily show that
$$\left(\widehat{Z}^{\textbf{z}_n}(t) - \frac{\textbf{z}_n}{\|\textbf{z}_n\|} - \frac{X(\|\textbf{z}_n\|t)}{\|\textbf{z}_n\|}, \ t\geq 0 \right) \overset{\mathcal{D}}{\Rightarrow} (\varphi^Z(t) - \varphi^Z(0) +  \nu t \textbf{1} , \ t\geq 0).$$

Now, given \eq{dynamic_equation} and \lemt{weak_convergence}, we have
$$(\varphi^\eta(t) B, \ t\geq 0) \overset{d}{=} (\varphi^Z(t) - \varphi^Z(0) +  \nu t \textbf{1} , \ t\geq 0).$$

By \rem{spike_rate}, the matrix $B$ is invertible and we have 
$$(\varphi^\eta(t), \ t\geq 0) \overset{d}{=} \left( \left( \varphi^Z(t) - \varphi^Z(0) +  \nu t \textbf{1}\right) B^{-1}, \ t\geq 0\right) .$$
Since $\varphi^\eta$ is a weak limit, we can assume without loss of generality $\varphi^\eta(t) = \left( \varphi^Z(t) - \varphi^Z(0) +  \nu t \textbf{1}\right) B^{-1}$. Thus, $\varphi^\eta$ is differentiable wherever $\varphi^Z$ is.

Assume that $\|\varphi^Z(t_0)\| > 0$ and $t_0$ is a regular point (see \sectn{fluid_model}). Let $N_0 = \sharp\{i: \ \varphi^Z_i(t_0) = 0\} < N$. Then, with a proper reordering, $ \varphi^Z_i(t_0) = 0$, for $i\in\{1, \ldots, N_0\}$ and $ \varphi^Z_i(t_0) > 0$, for $i\in\{N_0 + 1, \ldots, N\}$. Since $\varphi^{Z}_i(t) \geq 0$ and $t_0$ is a regular point, from $\varphi^{Z}_i(t_0) = 0$ we get $(\varphi^{Z}_i)'(t_0) = 0$. We find the values of
$$(\varphi^Z_i)'(t_0) = -\nu + H_i(\varphi^\eta_i)'(t_0) +  \sum_{j\neq i}w_{j}(\varphi^\eta_j)'(t_0).$$
We prove that $(\varphi^\eta_i)'(t_0) = 0$ for $i> N_0$ (if a potential is very far from the threshold then the neuron does not have a spike for a long time) and, therefore,
\begin{equation}\label{eq:zero_level_fluid}
0 = -\nu + (H_i-w_i)(\varphi^\eta_i)'(t_0)+ \sum_{j=1}^{N_0}w_{j}(\varphi^\eta_j)'(t_0), \ \ i=1,\ldots, N_0,
\end{equation}
\begin{equation}\label{eq:positive_level_fluid}
(\varphi^Z_i)'(t_0) = -\nu +  \sum_{j=1}^{N_0}w_j(\varphi^\eta_j)'(t_0), \ \ i=N_0+1, \ldots, N.
\end{equation}

Let 
$$h = \min_{N_0 + 1 \leq i \leq N} \varphi^Z_i(t_0).$$ 
We prove that for any $\Delta < h/(4\nu)$ and $i\in[N_0 +1, N]$ equality $\varphi^{\eta}_i(t_0 + \Delta) = \varphi^{\eta}_i(t_0)$ holds. Since $\widehat{Z}_i^{\textbf{z}_n}(t_0) \Rightarrow \varphi^Z_i (t_0)$, we have $\widehat{Z}_i^{\textbf{z}_n}(t_0) > h/2 > 2\nu\Delta$ a.s. for $n$ large. We have
\begin{equation*}
\begin{split}
\PB\{\eta^{\textbf{z}_n}_i(\|\textbf{z}_n\|(t_0 + \Delta)) > \eta^{\textbf{z}_n}_i(\|\textbf{z}_n\|t_0) \} & \leq \PB\{2\nu\Delta\|\textbf{z}_n\| + \inf_{0\leq s\leq \Delta} (X_i(\|\textbf{z}_n\|(t_0 + s)) - X_i(\|\textbf{z}_n\|t_0)) \leq 0\}\\
 & \leq \PB\{\nu\Delta\|\textbf{z}_n\| + \inf_{0\leq s\leq \Delta} (X^0_i(\|\textbf{z}_n\|(t_0 + s)) - X^0_i(\|\textbf{z}_n\|t_0)) \leq 0\}\\
 & = \PB\{\sup_{0\leq s\leq \Delta} X^0_i(\|\textbf{z}_n\|s) \geq \nu\Delta\|\textbf{z}_n\|\}.
\end{split}
\end{equation*}
Thus, by Theorem 36.8 from Sato (1999), we have convergence $\eta^{\textbf{z}_n}_i(\|\textbf{z}_n\|(t_0 + \Delta)) - \eta^{\textbf{z}_n}_i(\|\textbf{z}_n\|t_0) \to 0$ in probability and convergence $\widehat{\eta}^{\textbf{z}_n}_i(t_0 + \Delta) - \widehat{\eta}^{\textbf{z}_n}_i(t_0) \to 0$ a.s., as $n\to\infty$. Thus, equality $\varphi^{\eta}_i(t_0 + \Delta) = \varphi^{\eta}_i(t_0)$ holds for $\Delta < h/(4\nu)$ and $(\varphi^{\eta}_i)'(t_0) = 0$.

Using \rem{spike_rate} we solve system \eq{zero_level_fluid} and get
$$(\varphi^\eta_i)'(t_0) = \frac{\nu}{H_i - w_i} \frac{1}{1 + \sum_{k=1}^{N_0}\frac{w_k}{H_k - w_k}}, \ \ i=1,\ldots, N_0,$$ 
and therefore,
$$(\varphi^Z_i)'(t_0) = -\nu + \frac{\nu}{1 + \sum_{k=1}^{N_0}\frac{w_k}{H_k - w_k}} \sum_{j=1}^{N_0}\frac{w_j}{H_j - w_j} = -\frac{\nu}{1 + \sum_{k=1}^{N_0}\frac{w_k}{H_k - w_k}}, \ \ i=N_0+1, \ldots, N.$$

Therefore, the process $\varphi^Z$ is deterministic and piecewise linear. We have
$$\varphi_i^Z(0)\leq 1 \ \text{and} \ (\varphi^Z_i)'(t_0) \leq -\frac{\nu}{1 + \sum_{k=1}^{N}\frac{w_k}{H_k - w_k}}, \ \ i=1,\ldots, N,$$
for any regular point $t_0$ such that $\varphi_i^Z(t_0) >0$. Thus, from \eq{lipschitz_integration} we have that in time interval\\ $(0, \nu^{-1} (1 + \sum_{k=1}^{N}\frac{w_k}{H_k - w_k}))$ process $\varphi^Z$ reaches zero and stays there.

Let $\tau^{\textbf{z}}(\varepsilon, B) = \inf\{t> \varepsilon: \ Z^{\textbf{z}}(t) \in B\}$. Since fluid limits are stable, there exists $\kappa >0$ such that for $V = \{\textbf{z}\in \mathcal{Z}: \|\textbf{z}\| < k_0\}$ we have
$$\sup_{\textbf{z}\in V} \Expect\tau^{\textbf{z}}(\varepsilon, V) < \infty.$$

\section{Proof of \lemt{minorant_measure}}\label{sect:lemma_mixing}
We prove existence of a lower bound for $\inf_{\textbf{z} \in V} \PB\{Z^{\textbf{z}}(U) \in B\}$ where\\
 $V= \{\textbf{z}\in \mathcal{Z}: \|\textbf{z}\| < k_0\}$ (see the end of previous section).  By Theorem 19.2 from Sato (1999), the L\'evy process $X(t)$ can be represented as a sum $X^1(t) + X^2(t)$ of two independent processes, a jump process $X^1(t)$ and a Gaussian process $X^2(t)$ with drift. We consider cases where at least one coordinate is close enough to zero. If all the coordinates of $\textbf{z}$ are bounded away from zero then the proof follows similar lines.
 
Since random variables $\xi^{(1)}_{ij}$, $i, j\in [1, N]$, are strictly positive, there are constants $k^+_1, k^-_1 >0$ such that
$$p_1 = \PB\{A_1\} \equiv \PB\left\{(\xi^{(1)}_{ij})_{i, j=1}^N \in [k^-_1, k^+_1]^{N^2}\right\} >0.$$
Without loss of generality, we assume that $z_1 < k^-_1/6$.
 
First, we bound the jump process $X^1(t)$ in the time interval $[0, 2]$, which includes the time interval $[0, U]$, and take time instant $t_0 \leq 1/2$ such that $\PB\{X^1_1(t_0) < k^-_1/6\} >0$. Denote
$$A_2 =  \left\{\max_{1\leq i \leq N} (X^1_i(2)) < k_2\right\} \cap \left\{X^1_1(t_0) < k^-_1/6\right\},$$
and take a constant $k_2 >0$ such that $p_2  = \PB\{A_2\} >0.$ Next, we use the condition that the Gaussian process $X^2(t)$ is non-degenerate and none of its coordinates is a deterministic line. Thus, we denote
\begin{equation*}
\begin{split}
A_3 & = \left\{\max_{1\leq i \leq N} \sup_{0\leq t \leq \frac{1}{2}} (X^2_i(t)) \leq k_3-k_2 - k^-_1\right\} \cap \left\{\min_{1\leq i \leq N} \inf_{0\leq t \leq \frac{1}{2}} (X^2_i(t)) \geq -\frac{k^-_1}{2}\right\}\\
& \cap \left\{\inf_{0\leq t \leq t_0} (X^2_1(t)) \leq -\frac{k^-_1}{3}\right\}
\end{split}
\end{equation*}
and take a constant $k_3 > k_2 - k^-_1$ such that $p_3 = \PB\{A_3\} >0$. One can show that, given $A_2 \cap A_3$, the first spike occurs up to time $t_0$ and the second one can occur only after time $1/2$.

Denote the new set $D = A_1 \cap A_2 \cap A_3$. From independence of $X^1$, $X^2$ and $(\xi^{(1)}_{ij})_{i, j =1}^N$ we have $\PB\{D\} = p_1p_2p_3 >0$. We have
\begin{equation}\label{eq:away_from_zero_after_spike}
D \subseteq \left\{\frac{k^-_1}{2} \leq Z_i^{\textbf{z}}\left( \frac{1}{2}\right)  + X_i^1(U) - X_i^1\left( \frac{1}{2}\right) \leq k^+_1 + k_2 + k_3, \ i = 1, \ldots, N\right\}.
\end{equation}
We restrict ourselves to events without the second spike up to time $U$. Denote 
	$$G_2(t_1, t_2, k) = \left\{\min_{1\leq i \leq N} \inf_{t_1\leq s \leq t_2} (X^2_i(s)) > -k\right\}.$$ 
Using \eq{away_from_zero_after_spike}, we get that, given $G_2\left( 1/2, U, k^-_1/2\right) \cap D$, the second spike occurs after time $U$. 

Let $K = k^+_1 + k_2 + k_3$. We prove that, for any point $\textbf{y} \in (K, 2K)^N$ and a measurable set $\Delta \subset [0, K]^N$, there is a number $p>0$ such that
$$
\PB\left\{\left\{Z^{\textbf{z}}(U) \in \textbf{y} + \Delta\right\} \cap G_2\left( \frac{1}{2}, U, \frac{k^-_1}{2}\right) \cap D \right\}  \\
\geq p \lambda(\Delta),
$$
where $\lambda$ is the Lebesgue measure. Denote, $\widehat{y}_i = y_i - Z_i^{\textbf{z}}\left( 1/2\right)  - \left( X_i^1(U) - X_i^1\left( 1/2\right)\right) $, for $i\in[1, N]$. Then we have
$$\left\{Z^{\textbf{z}}(U) \in \textbf{y} + \Delta\right\} = \left\{ X^2(U) - X^2\left( \frac{1}{2} \right)  \in\widehat{\textbf{y}} + \Delta\right\}.$$

Since $X^2$ is a Markov process, the events $ \left\{ X^2(U) - X^2\left( 1/2 \right)  \in\widehat{\textbf{y}} + \Delta\right\} \cap G_2\left( 1/2, U, k^-_1/2\right)$ and $D$ are independent, conditioned on a value of $\widehat{\textbf{y}}$. Thus, we have
\begin{equation*}
\begin{split}
&\PB\left\{\left\{ X^2(U) - X^2\left( \frac{1}{2} \right)  \in\widehat{\textbf{y}} + \Delta\right\} \cap G_2\left( \frac{1}{2}, U, \frac{k^-_1}{2}\right) \cap D \right\}\\
& = \Expect\left( \PB\left\{\left\{ X^2(U) - X^2\left( \frac{1}{2} \right)  \in\widehat{\textbf{y}} + \Delta\right\} \cap G_2\left( \frac{1}{2}, U, \frac{k^-_1}{2}\right) \mid \ \widehat{\textbf{y}}\right\} \PB\left\{D \mid \ \widehat{\textbf{y}}\right\}\right) 
\end{split}
\end{equation*}

Next, we need a technical lemma regarding a monotonicity property of the Brownian bridge.
\begin{Lemma}\label{lem:brownian_bridge}
For any $t, k>0$ and $\Delta \subset [0, \infty)^N$ we have
$$\PB\{G_2(0, t, k)  \mid \ X^2(t) \in \Delta\} \geq \PB\{G_2(0, t, k) \mid \ X^2(t) = \textbf{0}\} >0.$$
\end{Lemma}

Given $D$, we have $\widehat{y}_i \geq 0$, for $i \in [1, N]$, and we can use \lemt{brownian_bridge} to obtain
\begin{multline*}
\PB\left\{\left\{ X^2(U) - X^2\left( \frac{1}{2} \right)  \in\widehat{\textbf{y}} + \Delta\right\} \cap G_2\left( \frac{1}{2}, U, \frac{k^-_1}{2}\right) \mid \ \widehat{\textbf{y}}\right\}\\
 \overset{a.s.}{\ge} \PB\left\{ X^2(U) - X^2\left( \frac{1}{2} \right)  \in\widehat{\textbf{y}} + \Delta\mid \ \widehat{\textbf{y}}\right\} \PB\left\lbrace G_2\left( \frac{1}{2}, U, \frac{k^-_1}{2}\right) \mid \ X^2\left( U\right) -X^2\left( \frac{1}{2}\right) = \textbf{0}  \right\rbrace. 
\end{multline*}

The density of $X^2(t)$  is bounded away from zero on any compact set, and $X^2$ and $U$ are independent. Therefore, there exists $p_4 >0$ such that, given $D$, for a measurable set $\Delta \subseteq [0, K]^N$  we have
$$\PB\left\{ X^2(U) - X^2\left( \frac{1}{2} \right)  \in\widehat{\textbf{y}} + \Delta\mid \ \widehat{\textbf{y}}\right\} \overset{a.s.}{\geq} p_4\lambda(\Delta).$$
Denote $p'_4 = p_4\PB\left\lbrace G_2\left( \frac{1}{2}, U, \frac{k^-_1}{2}\right) \mid \ X^2\left( U\right) -X^2\left( \frac{1}{2}\right) = \textbf{0}  \right\rbrace >0$. Combining altogether, we get that if $z_1 < k^-_1/6$ then
$$\PB\{Z^{\textbf{z}}(U) \in \textbf{y} + \Delta\} \geq p_1p_2p_3p'_4\lambda(\Delta),$$
for $\textbf{y} \in (K, 2K]$ and $\Delta \subseteq [0, K]^N$.

\appendix
\section*{Appendix}
\section*{Comments on \rem{partial_stability}}
In a system of two inhibitory neurons it is sufficient for stability to assume that the signals are smaller on the average than the thresholds. However, in a system of three inhibitory neurons it is not enough. For the matrix
$B=\begin{pmatrix}
8 & 2 & 6\\
2 & 8 & 6\\
6 & 6 & 8
\end{pmatrix}$  and drifts $\nu_i = -1$, $i=1,2,3,$, first two neurons can form a stable system that 'pushes' the potential of the third neuron to infinity. Here is an example of sufficient conditions on matrix $B$ and parameters $\nu_i$ to avoid such cases (we believe that they can be weakened).
\begin{itemize}
\item For every set $S \subseteq [1, N]$ the matrix $B^S = (b_{ij})_{i, j \in S}$ is invertible and $a^S_i = \left( (B^S)^{-1}f^S\right)_i >0 $, $i\in [1, N]$ where $f^S = (\nu_i)_{i\in S}$;
\item $\sum_{i\in S} a^S_i \sum_{j\notin S} b_{ij} < \sum_{j\notin S} \nu_j.$
\end{itemize}

\section*{Proof of \rem{spike_rate}}

We prove the remark by assuming that the system of equations $\textbf{x} B = \textbf{1}$ has a solution, and then we prove that the solution is unique and has all positive coordinates. Let us rewrite the system  $\textbf{x} B = \textbf{1}$ as
$$\sum_{j=1}^N x_j  b_{ji} = H_i x_i + \sum_{j\neq i}  w_j x_j = (H_i-w_i)x_i + \sum_{j=1}^N w_j x_j = 1 \ \ i= 1, \ldots, N.$$

Denote $M = \sum_{j=1}^N w_j x_j$ and get
$$x_i = \frac{1 - M}{H_i - w_i} \ \ i= 1, \ldots, N \ \Rightarrow \ M = (1-M)\sum_{i=1}^N\frac{w_i}{H_i - w_i} \ \Rightarrow \ M = \frac{\sum_{i=1}^N\frac{w_i}{H_i - w_i}}{1 + \sum_{i=1}^N\frac{w_i}{H_i - w_i}} < 1.$$
Thus, $M$ and, therefore, $\textbf{x}$ are uniquely defined through $\{H_i, w_i\}_{i=1}^n$ and $x_i > 0$, $i=1, \ldots, N$. The rest of the proof follows from the proof of \lemt{positive_recurrence}.

\section*{Proof of \lemt{brownian_bridge}}

First, take $X(t)$, $t\in [0,1]$, a standard one-dimensional Brownian motion. Then for the Brownian bridge $B_x(t) = xt + B_0(t)$ and $x\leq y$ we have
\begin{equation*}
\begin{split}
\PB\{ \inf_{0\leq t \leq 1} (X(t)) \geq -k | \ X(1) = x\} &  = \PB\{ \inf_{0\leq t \leq 1} (xt + B_0(t)) \geq -k\}\\
 & \leq \PB\{ \inf_{0\leq t \leq 1} (yt + B_0(t)) \geq -k\}\\
& = \PB\{ \inf_{0\leq t \leq 1} (X(t)) \geq -k | \ X(1) = y\}.
\end{split}
\end{equation*}
For a general $N$-dimensional Brownian motion $X(t)$, $t\in[0, 1]$, with a non-singular covariance matrix $\Sigma$, there exists an invertible matrix $L$ and a vector $\textbf{v}$ such that $W(t) = LX(t) + \textbf{v}t$ is a vector of $N$ independent standard Brownian motions. Let $B_{\textbf{0}}(t)$ denote the corresponding $N$-dimensional Brownian bridge. Denote $A_k = [-k, \infty)^N$. Then for $\textbf{x}, \textbf{y}$, such that $x_i \leq y_i$, we have
\begin{equation*}
\begin{split}
&\PB\left\{\min_{1 \leq i \leq N} \inf_{0\leq t \leq 1} (X_i(t)) \geq -k | \ X(1) = \textbf{x} \right\}\\
& = \PB\left\{X(t) \in A_k, \  t\in[0, 1]| \ X(1) = \textbf{x}\right\} = \PB\left\{W(t) \in LA_k + \textbf{v}t, \ t\in[0, 1] | \ W(1) = L\textbf{x} + \textbf{v}\right\}\\
& = \PB\left\{L\textbf{x}t + \textbf{v}t + B_{\textbf{0}}(t) \in LA_k + \textbf{v}t, \ t\in[0, 1]\right\}
= \PB\left\{\textbf{x}t + L^{-1}B_{\textbf{0}}(t) \in [-k, \infty)^N, \ t\in[0, 1]\right\}\\
& \leq \PB\left\{\textbf{y}t + L^{-1}B_{\textbf{0}}(t) \in [-k, \infty)^N, \ t\in[0, 1]\right\} = \PB\left\{\min_{1 \leq i \leq N} \inf_{0\leq t \leq 1} (X_i(t)) \geq -k | \ X(1) = \textbf{y}\right\}.
\end{split}
\end{equation*}
Using the properties of the Brownian bridge, one can replace the $\textbf{y}$ with a measurable set $\Delta$ and prove the statement.

\section*{Acknowledgements}

I would like to thank Thibauld Taillefumier for the introduction to the model and overview of the problem. I would also like to thank my supervisors Sergey Foss and Seva Shneer for insightful remarks and ideas.

\end{document}